\documentclass[11pt]{article}
\usepackage{amsmath}
\usepackage{amsfonts}
\usepackage[bookmarksnumbered=true]{hyperref}

\newtheorem{theorem}{Theorem}[section]
\newtheorem{proposition}[theorem]{Proposition}
\newtheorem{lemma}[theorem]{Lemma}
\newtheorem{definition}[theorem]{Definition}
\newtheorem{remark}[theorem]{Remark}
\newtheorem{example}[theorem]{Example}
\def\dint{\displaystyle\int}

\begin{document}

\title{Multidimensional BSDEs with uniformly continuous coefficients: the general result
\thanks{Partially supported by the National Natural Science Foundation of China (Grant Nos. 11371226, 11071145, 11301298 and 11231005), the Foundation for Innovative Research Groups of National Natural Science Foundation of China (Grant No. 11221061) and the 111 Project (Grant No. B12023).}}
\author{Min Li$^{\rm a}$\thanks{E-mail: minli2006@163.com} and Yufeng Shi$^{\rm a, \rm b}$\thanks{Corresponding author. E-mail: yfshi@sdu.edu.cn}\\
{\small$^{\rm a}$Institute for Financial Studies, Shandong University, Jinan 250100, China}\\
{\small$^{\rm b}$School of Statistics, Shandong University of Finance and Economics, Jinan 250014, China}}
\maketitle

\begin{abstract}In this paper, by introducing a new notion of \textit{envelope} of the stochastic process, we construct a family of random differential equations whose solutions can be viewed as solutions of a family of ordinary differential equations and prove that the multidimensional backward stochastic differential equations (BSDEs for short) with the general uniformly continuous coefficients are uniquely solvable. As a result, we solve the open problem of multidimensional BSDEs with uniformly continuous coefficients.\\

\indent{\it keywords:} Backward stochastic differential equations, envelope, Girsanov's theorem.\\


\end{abstract}

\section{Introduction}\label{sec:1}

In this paper we study the multidimensional backward stochastic differential equations (BSDEs for short) of the following form
\begin{equation}\label{eq:1.1}
y_t=\xi+\int_t^1f(s,y_s,z_s){\rm d}s-\int_t^1z_s{\rm d}B_s,\quad 0\leq t\leq 1,
\end{equation}
where $\xi$ is an $R^d$-valued random variable and $(B_t)_{0\leq t\leq 1}$ is an $m$-dimensional standard Brownian motion. If $\xi$ is square integrable and the driver $f$ is Lipschitz continuous in $(y,z)$,  Pardoux and Peng \cite{r12} proved the existence and uniqueness of the solution to BSDE (\ref{eq:1.1}). Since then, many efforts have been made to relax the Lipschitz conditions on the generator $f$. In one-dimensional case, many results for BSDEs with more general generators have been obtained with the help of the comparison theorems which can be referred to El Karoui et al. \cite{r2}. For example, Lepeltier and San Martin \cite{r7,r8} solved the BSDEs with continuous coefficients and superlinear quadratic coefficients. Kobylanski \cite{r6} proved the solvability of BSDEs with quadratic growth coefficients, and so on.

For the multidimensional situation, the solvability of BSDE (\ref{eq:1.1}) with non-Lipschitzian generators becomes complicated due to that there is no multidimensional result directly analogous to the one-dimensional comparison theorem. Nevertheless, when the driver $f$ is Lipschitz with respect to $z$ and non-Lipschitzian in $y$, Fan et al. \cite{r4}, Mao \cite{r10} and Pardoux \cite{r11} proved the existence and uniqueness of the solutions. And for the case that $f$ in $z$ is also non-Lipschitzian, precisely speaking, under the following uniformly continuous conditions:\\
(i) $(y,z)\mapsto f(t,y,z)$ is uniformly continuous uniformly in $(\omega,t)$ and satisfies Assumption 2 below;\\
(ii) the $i$th component $f_i$ of $f$ depends only on the $i$th row of $z$,\\
Hamad\`{e}ne \cite{r5} obtained the existence of a solution and the uniqueness result was given by Fan et al. \cite{r3}.

In this paper we will show that BSDE (\ref{eq:1.1}) also has a unique solution without the condition (ii). Actually, it is not easy. Referring to Hamad\`{e}ne \cite{r5} and Fan et al. \cite{r3}, we can see that the Girsanov theorem, which is guaranteed by the condition (ii), plays an important role in proving the existence and uniqueness of the solution. In general, without the condition (ii), one is not able to use the Girsanov theorem on the whole interval $[0,1]$. In this paper, we use a new argument to successfully overcome these difficulties and accordingly give a more complete solution to the open problem of multidimensional BSDEs with uniformly continuous coefficients.

Thanks to the condition (i), we can introduce a sequence of Lipschitz functions $(f_n)_{n\geq 0}$ which converges to $f$ uniformly on the whole space of $(y,z)$. As usual, let $(y_n,z_n)_{n\geq 0}$ be the solutions to BSDE (\ref{eq:1.1}) associated with $(f_n,\xi)_{n\geq 0}$ and we will show that $(y_n,z_n)_{n\geq 0}$ converge to the solution of BSDE (\ref{eq:1.1}). In order to prove the convergence of $(y_n,z_n)_{n\geq 0}$, for all $(y_m,z_m)$ and $(y_n,z_n)$, $m\neq n$, we will construct countable stopping time intervals. For any given stopping time interval which is already constructed, we further construct a stopping time subinterval on it, which enables us to apply Girsanov's theorem. Then with the help of the envelope of the stochastic process introduced in the following, we construct a family of random backward differential equations, which can be viewed as BSDEs and backward ordinary differential equations respectively, such that their solutions can dominate $|y_n-y_m|$ on the constructed stopping time subinterval. By changing the time horizon and the terminal values of the random backward differential equations, we show that the sequence $(y_n,z_n)_{n\geq 0}$ will converge to the solution of BSDE (\ref{eq:1.1}). For the uniqueness of the solution, we take the similar procedure.

The paper is organized as follows. In Section 2 we introduce some notations and give some preliminary results. In Section 3 we prove our main results.

\section{Notations and Preliminary}\label{sec:2}

Let $(\Omega,{\cal{F}},P)$ be a complete probability space and $\{B_t,0\leq t\leq 1\}$ be an $m$-dimensional standard Brownian motion defined on $(\Omega,{\cal{F}},P)$. We denote by $F=\{{\cal{F}}_t,0\leq t\leq 1\}$ the natural filtration of $(B_t)_{0\leq t\leq 1}$, where ${\cal{F}}_0$ contains all $P$-null sets of $\cal{F}$. The terminal condition $\xi$ is an ${\cal{F}}_1$-measurable $d$-dimensional random vector and the driver $f:[0,1]\times\Omega\times R^d\times R^{d\times m}\rightarrow R^d$ is a ${\cal{P}}\otimes {\cal{B}}(R^d)\otimes {\cal{B}}(R^{d\times m})$-measurable mapping, where $\cal{P}$ denotes the $\sigma$-algebra of ${\cal{F}}_t$-progressive measurable sets on $[0,1]\times\Omega$. In this paper, we need the following notation:\\

$\bullet$ $H^{2,k}$: all $\cal{P}$-measurable processes $U={(U_t)}_{0\leq t\leq 1}$ with values in $R^k$ such that $E[\int_0^1|U_s|^2ds]<+\infty$;\\

$\bullet$ $H^{2,k\times l}$: all $\cal{P}$-measurable processes $U={(U_t)}_{0\leq t\leq 1}$ with values in $R^{k\times l}$ such that $E[\int_0^1|U_s|^2ds]<+\infty$;\\

$\bullet$ $S^{2,k}$: all continuous $\cal{P}$-measurable processes $U={(U_t)}_{0\leq t\leq 1}$ with values in $R^k$ such that $E[\sup\limits_{0\leq t\leq 1}|U_t|^2]<+\infty$.\\

\textbf{Assumption 1}.\quad  The process $f(t,\omega,0,0)$, $(0\leq t\leq 1)$ belongs to $H^{2,d}$ and, for any $(y,z)\in R^d\times R^{d\times m}$, the process $(f(t,\omega,y,z))_{0\leq t\leq 1}$ is $\cal{P}$-measurable.\\

\textbf{Assumption 2}.\\
(i) $f$ is uniformly continuous in $y$ uniformly with respect to $(t,\omega,z)$, i.e., there exists a continuous non-decreasing function $\Phi$ from $R^+$ to $R^+$ with at most linear growth and satisfying $\Phi(0)=0$ and $\Phi(x)>0$ for all $x>0$ such that:
\[
|f(t,y_1,z)-f(t,y_2,z)|\leq \Phi(|y_1-y_2|), \quad a.s.,\quad \forall t,y_1,y_2,z,
\]
moreover, $\int_{0^+}[\Phi(x)]^{-1}dx=+\infty$.\\
(ii) $f$ is uniformly continuous in $z$, i.e., there exists a continuous function $\Psi$ from $R^+$ into itself with at most linear growth and satisfying $\Psi(0)=0$, such that:
\[
|f(t,y,z_1)-f(t,y,z_2)|\leq \Psi(||z_1-z_2||),\quad a.s.,\quad \forall t,y,z_1,z_2,
\]
where $||z||=[tr(zz^*)]^{\frac{1}{2}}$, $z^*$ is the transpose of $z$.\\

\begin{definition}
\label{de:2.1}
A pair of $\cal{P}$-measurable process $(y,z)=(y_t,z_t)_{0\leq t\leq 1}$ valued in $R^d\times R^{d\times m}$ is called an adapted solution of BSDE (\ref{eq:1.1}), if $(y,z)\in S^{2,d}\times H^{2,d\times m}$ and satisfies BSDE (\ref{eq:1.1}).
\end{definition}

\begin{proposition}
\label{pro:2.1}
(Hamad\`{e}ne \cite{r5}) Let $\phi$ be a continuous non-decreasing function from $R^+$ to $R^+$ such that $\int_{0^+}\frac{dx}{\phi(x)}=+\infty$ and $\phi(x)\leq ax+b$ for all $x\in R^+$, where $a$ and $b$ are given non-negative constants. If $\phi(x)>0$ for all $x>0$, then the following backward differential equation
\[
u^{\gamma,\varepsilon}(t)=\gamma+\int_t^1(\phi(u^{\gamma,\varepsilon}(s))+\varepsilon)ds,\quad 0\leq t\leq 1,
\]
has a unique solution, where $\gamma\geq 0,\varepsilon\geq 0$. In particular, the solution $u^{\gamma,\varepsilon}(t)$ is continuous with respect to $\gamma$, $\varepsilon$ and satisfies that $\lim\limits_{\varepsilon\rightarrow 0}\lim\limits_{\gamma\rightarrow 0}u^{\gamma,\varepsilon}(t)=0$.
\end{proposition}

Now, we introduce a notion of envelope which plays an important role in this paper. By Lemma 3 in Lepeltier and San Martin \cite{r8}, we know, for a given constant $a\in R^+$, the following BSDE
\begin{equation}\label{eq:2.1}
y(t)=a+\int_t^1(\phi(y(s))+\varepsilon)ds-\int_t^1z(s)dB(s),
\end{equation}
has a unique solution $y\equiv x$ and $z\equiv 0$, where
\[
x(t)=a+\int_t^1(\phi(x(s))+\varepsilon)ds.
\]
For all stopping time $\tau$ such that $0\leq \tau\leq 1$, it follows that
\begin{equation*}
\begin{array}{lll}
y(t)&=&y(\tau)+\dint_t^{\tau}(\phi(y(s))+\varepsilon)ds-\dint_t^{\tau}z(s)dB(s)\\
&=&x(\tau)+\dint_t^{\tau}(\phi(y(s))+\varepsilon)ds-\dint_t^{\tau}z(s)dB(s),\quad 0\leq t\leq\tau,
\end{array}
\end{equation*}
and
\begin{equation}\label{eq:2.2}
x(t)=x(\tau)+\int_t^{\tau}(\phi(x(s))+\varepsilon)ds,\quad 0\leq t\leq\tau.
\end{equation}
Obviously, it still holds that $x(t)=y(t), 0\leq t\leq\tau,$ a.s.. On the other hand, Eq. (\ref{eq:2.2}) can also be viewed as a random backward differential equation on $[0,\tau]$, and by Proposition \ref{pro:2.1}, for P-a.s. $\omega\in\Omega$, Eq. (\ref{eq:2.2}) has a unique solution $x(t,\omega), 0\leq t\leq\tau(\omega)$.

Let $X(t), 0\leq t\leq 1,$ be a nonnegative continuous and bounded process such that $\sup\limits_{0\leq t\leq 1}X(t)\leq C$, a.s., where $C$ is a nonnegative constant, and let $u^{\gamma}(t)$ be the solution of the following backward differential equation:
\begin{equation}\label{eq:2.3}
u^{\gamma}(t)=\gamma+\dint_t^1(\phi(u^{\gamma}(s))+\varepsilon)ds.
\end{equation}
Now, we denote $\phi_{\varepsilon}(x):=\phi(x)+\varepsilon$ and, for a given stopping time $\tau$ such that $0\leq \tau \leq 1$,
\begin{equation}\label{eq:2.4}
\begin{array}{rr}
\Upsilon:=\left\{\gamma\geq 0:u^{\gamma}(\tau)\geq X(\tau), u^{\gamma} \text{ is the }\right.\\
\left. \text{ solution of Eq. (\ref{eq:2.3})}\right\}.
\end{array}
\end{equation}
If $\gamma\geq C$, it is easy to obtain that $X(\tau)\leq C\leq u^{\gamma}(\tau), a.s.$. This implies that $\Upsilon$ is a nonempty set. Denoting $\gamma^0:=\text{inf}\{\gamma:\gamma\in\Upsilon\}$, then by the continuity of the process $X$ and the continuity of $u^{\gamma}(t)$ with respect to $\gamma$, we have $u^{\gamma^0}(\tau)\geq X(\tau)$, it follows that $\gamma^0\in\Upsilon$. Furthermore, we denote
\begin{equation}\label{eq:2.5}
\begin{array}{rr}
\theta(\tau):=\text{essinf}\left\{u^{\gamma}(\tau): u^{\gamma}(\tau)\geq X(\tau), u^{\gamma} \text{ is the }\right.\\
\left. \text{ solution of Eq. (\ref{eq:2.3})}, \gamma\geq 0\right\},
\end{array}
\end{equation}
then\\

\begin{definition}
\label{def:2.2}
For a given stopping time $\tau$ such that $0\leq\tau\leq 1$, $\theta(\tau)$ denoted by (\ref{eq:2.5}) is called the $\phi_{\varepsilon}$-envelope of the process $X$ at $\tau$.
\end{definition}

\begin{remark}
\label{rmk:2.1}
Obviously, we have $\theta(\tau)=u^{\gamma_0}(\tau)$. By the continuity of the process $X$ and the continuity of $u^{\gamma}(t)$ with respect to $\gamma$, it follows that, for any stopping times $(\tau_i)_{i\in I}$, where $I$ denotes an index set, such that $\tau_i\rightarrow \tau$, $\theta(\tau_i)\rightarrow\theta(\tau)$. In particular, if $X(1)=0$, we have $\theta(1)=0$.
\end{remark}

\begin{remark}
\label{rmk:2.2}
For all $\alpha>0$, if we denote $X_{\alpha}(t):=\alpha X(t)$ and the $\phi_{\varepsilon}$-envelope of the process $X_{\alpha}$ at $\tau$ by $\theta_{\alpha}(\tau)$, then we have $\theta_{\alpha}(\tau)=\alpha\theta(\tau)$. And, as $\alpha\rightarrow 1$, it follows that $\theta_{\alpha}(\tau)\rightarrow \theta(\tau)$.
\end{remark}

\begin{lemma}
\label{lem:2.1}
Let ${\cal{A}}=\{\nu_i\}_{i\in I}$ be a set of stopping times such that $\nu_{i_1}\wedge\nu_{i_2}\in {\cal{A}}, \forall i_1, i_2\in I$. If we denote $\tau:=\text{essinf}_{i\in I}\nu_i$, then $\tau$ is a stopping time.
\end{lemma}
{\bf Proof.}
By the fact that $\nu_{i_1}\wedge\nu_{i_2}\in {\cal{A}}, \forall i_1, i_2\in I$, we can obtain that there exists a sequence of decreasing stopping times $(\tau_n)_{n\geq 1}\subset {\cal{A}}$ such that
\[
\lim\limits_{n\rightarrow +\infty}E[\tau_n]=\inf\limits_{i\in I}E[\nu_i].
\]
Denoting $\eta:=\wedge_{n\in N}\tau_n$ and $\eta_n:=\tau_n\wedge\nu, \nu\in {\cal{A}}$, then we have $\lim\limits_{n\rightarrow +\infty}\eta_n=\eta\wedge\nu$, and
\[
E[\eta\wedge\nu]=\lim\limits_{n\rightarrow\infty}E[\eta_n]\geq\inf\limits_{i\in I}E[\nu_i]=E[\eta].
\]
From the fact that $\eta\wedge\nu\leq \eta$, we have $\eta\wedge\nu=\eta$, a.s., that is P-a.s. $\eta\leq \nu, $. On the other hand, for any random variable $\theta$ satisfying $\theta\leq \nu_i, \forall i\in I$, it is obviously $\eta\geq \theta$, a.s., hence, we have $\tau=\eta$, a.s., and the desired result is obtained.
\hfill $\Box$

\begin{proposition}
\label{pro:2.2}
(Hamad\`{e}ne \cite{r5}) Let $g$ be a mapping $(t,\omega,y,z)\mapsto g(t,\omega,y,z)$ from $[0,1]\times\Omega\times R^d\times R^{d\times m}$ to $R^d$ satisfying Assumption 1. In addition, there exists a continuous function $\varphi$ from $R^+\times R^+$ into $R^+$ such that $\varphi(0,0)=0$ and
\[
|g(t,\omega,y_1,z_1)-g(t,\omega,y_2,z_2)|\leq \varphi(|y_1-y_2|,||z_1-z_2||), \quad\forall t,y_1,y_2,z_1,z_2, a.s.
\]
Then there exists a sequence $(g_n)_{n\geq 0}$ such that:\\
(i) For any $n\geq 0$, $g_n$ is a mapping from $[0,1]\times\Omega\times R^d\times R^{d\times m}$ into $R^d$ satisfying Assumption 1 and which is Lipschitz with respect to $(y,z)$ uniformly in $(t,\omega)$.\\
(ii) For all $\varepsilon>0$, there is an $N_{\varepsilon}\geq 0$ such that, $\forall n\geq N_{\varepsilon}$, $|g_n(t,\omega,y,z)-g(t,\omega,y,z)|\leq \varepsilon$ for all $t,y,z$, a.s.
\end{proposition}

\begin{remark}
\label{rmk:2.3}
Actually, from Hamad\`{e}ne \cite{r5}, we have $g_n:=g\ast\psi_n$, the convolution product of $g$ and $\psi_n$, where $\psi_n:(y,z)\in R^{d+d\times m}\mapsto\psi_n(y,z)=n^2\psi(ny,nz)$ and $\psi$ is a function of $C^{\infty}(R^{d+d\times m},R^+)$ with compact support and satisfies $\int_{R^{d+d\times m}}\psi(y,z)dydz=1$.
\end{remark}

\section{The main result}

For a given process $h\in H^{2,1}$, if we denote $H(t):=\int_0^t|h(s)|ds$, we have

\begin{lemma}
\label{lem:3.1}
There exist countable stopping times $\{\tau_k^r\}_{r\in N,k\in N}$ such that, for P-a.s.  $\omega\in\{\omega:\tau_k^r(\omega)<\tau_{k+1}^r(\omega)\}$, $H(\omega,\tau_k^r(\omega))=H(\omega,\tau_{k+1}^r(\omega))$ or $H(\omega,\cdot)$ is strictly increasing on $[\tau_k^r(\omega),\tau_{k+1}^r(\omega)]$. Moreover, for P-a.s. $\omega\in \Omega$,  $\bigcup\limits_{r\in N }\bigcup\limits_{k\in N }[\tau_k^r(\omega),\tau_{k+1}^r(\omega)]=[0,1], a.e.$
\end{lemma}

Before the proof of Lemma \ref{lem:3.1}, we first see the following example:

\begin{example}
\label{ex:3.1}
Let $\{h_n\}_{n\geq 1}$ be a sequence of ${\cal{P}}$-measurable processes defined on $[0,1]$ satisfying $h_{2n}\equiv 0$ and
\begin{equation*}
h_{2n-1}(t)=
\left\{
\begin{array}{ll}
|B_t|+1,\quad t\in[\frac{1}{2^{2n-1}},\frac{1}{2^{2n-2}}],\\
0,\qquad t\in[0,\frac{1}{2^{2n-1}})\cup(\frac{1}{2^{2n-2}},1].
\end{array}
\right.
\end{equation*}
If we denote $h(t):=\sum_{n=1}^{\infty}h_n(t)$, $H(t):=\int_0^th(s)ds$ and
\begin{equation*}
\begin{array}{lll}
\pi_0:=0,\\
\pi_{2l-1}(\omega):=\inf\{t>\pi_{2l-2}(\omega):H(t,\omega)=H(s,\omega), \exists s\in (t,1]\}\wedge 1,\\
\pi_{2l}(\omega):=\inf\{t>\pi_{2l-1}(\omega):H(t,\omega)<H(s,\omega), \forall s\in (t,1]\}\wedge 1,\quad l\in N,
\end{array}
\end{equation*}
we have P-a.s. $\pi_l=0$, $\forall l\in N$.
\end{example}

\begin{remark}
\label{rek:3.1}
Actually, for all stopping times $\theta_1, \theta_2$ satisfying $\theta_1\leq\theta_2$ and $P(\theta_1<\theta_2)>0$, there exists a process on $[\theta_1,\theta_2]$ satisfying the same property of $h$ in Example \ref{ex:3.1}.
\end{remark}

Now, we give the proof of Lemma \ref{lem:3.1}.\\

{\bf Proof.}¡¡¡¡¡¡
Without loss of generality, we can suppose $P(\{\omega:H(\omega)\equiv 0\})<1$. Now, we denote $\tau_r(\omega):=\inf\{s:H(\omega,s)\geq r\}\wedge 1, r\in Q^+$, where $Q^+$ denotes the nonnegative rational numbers set on $R$, and
\begin{equation}\label{eq:3.1}
\begin{array}{lll}
\pi_0^r:=\tau_r,\\
\pi_{2l-1}^r(\omega):=\inf\{t>\pi_{2l-2}^r(\omega):H(\omega,t)=H(\omega,s), \exists s\in (t,1]\}\wedge 1,\\
\pi_{2l}^r(\omega):=\inf\{t>\pi_{2l-1}^r(\omega):H(\omega,t)<H(\omega,s), \forall s\in (t,1]\}\wedge 1, \quad l\in N.
\end{array}
\end{equation}
Obviously, $\tau_r$ and $\pi_l^r, l\in N,$ are all stopping times. If for arbitrary $r\in Q^+$, it holds that $\pi_l^r=\tau_r, \forall l\in N$, then, from the continuity of $H$, we can obtain $H\equiv 0$, a.s.. This contradicts to the fact that $P(\{\omega:H(\omega)\equiv 0\})<1$. Hence, there must exist $r_0\in Q^+$ such that $P(\pi_1^{r_0}>\tau_{r_0})>0$. From the definition of $\{\pi_l^r\}_{l\geq 1}$, we have, for P-a.s. $\omega\in\{\omega:\pi_{2l-1}^r(\omega)<\pi_{2l}^r(\omega)\}$, $H(\omega,\pi_{2l-1}^r(\omega))=H(\omega,\pi_{2l}^r(\omega))$ and, for P-a.s. $\omega\in\{\omega:\pi_{2l-2}^r(\omega)<\pi_{2l-1}^r(\omega)\}$, $H(\omega,\cdot)$ is strictly increasing on $[\pi_{2l-2}^r(\omega),\pi_{2l-1}^r(\omega)]$. Then, from the continuity of $H$ and the fact that $Q^+$ is countable, we get the desired result.
\hfill $\Box$

\begin{remark}
\label{rmk:3.3}
For a given $r_0\in Q^+$, the set of stopping times $\{\pi_l^{r_0}\}_{l\geq 1}$ defined as (\ref{eq:3.1}) may be finite. Namely, there exists a constant $l_0\geq 1$ such that $\pi_l^{r_0}=\pi_{l_0}^{r_0}, \forall l>l_0$. For example, we consider the process $H$ in Example \ref{ex:3.1} and denote $\tau_{r_n}(\omega):=\inf\{s:H(s)\geq r_n\}\wedge 1, r_n\in \{r_n\}_{n\geq 1}$, where $\{r_n\}_{n\geq 1}$ is a subset of $Q^+$ satisfying $r_n>0, \forall n,$ and $r_n\rightarrow 0$ as $n\rightarrow +\infty$. Now, we denote
\begin{equation*}
\begin{array}{lll}
\pi_0^{r_n}:=\tau_{r_n},\\
\pi_{2l-1}^{r_n}(\omega):=\inf\{t>\pi_{2l-2}^{r_n}(\omega):H(\omega,t)=H(\omega,s), \exists s\in (t,1]\}\wedge 1,\\
\pi_{2l}^{r_n}(\omega):=\inf\{t>\pi_{2l-1}^{r_n}(\omega):H(\omega,t)<H(\omega,s), \forall s\in (t,1]\}\wedge 1,\quad l\in N.
\end{array}
\end{equation*}
Obviously, $\forall r_n>0$, the set $\{\pi_l^{r_n}\}_{l\geq 1}$ is finite.
\end{remark}

\begin{remark}
\label{rmk:3.4}
If, for P-a.s. $\omega\in\{\omega:\tau_k^r(\omega)<\tau_{k+1}^r(\omega)\}$, $H(\omega,\tau_k^r(\omega))=H(\omega,\tau_{k+1}^r(\omega))$, then, on $[\tau_k^r,\tau_{k+1}^r]$, we can take a version of the process $h$ such that $h=0$, a.s. And for the case of $H(\omega,\cdot)$ is strictly increasing on $[\tau_k^r(\omega),\tau_{k+1}^r(\omega)]$, we can also take a version of $h$ on $[\tau_k^r,\tau_{k+1}^r]$ such that $|h|>0$, a.s. Therefore, there exists a version of $h$ on $[0,1]$ such that, for P-a.s. $\omega\in\{\omega:\tau_k^r(\omega)<\tau_{k+1}^r(\omega)\}$, $h(\omega,\cdot)=0$ or $|h(\omega,\cdot)|>0$ on $[\tau_k^r(\omega),\tau_{k+1}^r(\omega)]$. At this time, the stopping times defined by (\ref{eq:3.1}) can be denoted by
\begin{equation*}
\begin{array}{lll}
\pi_{2l-1}^r(\omega):=\inf\{t>\pi_{2l-2}^r(\omega):|h(\omega,t)|=0\}\wedge 1, \\
\pi_{2l}^r(\omega):=\inf\{t>\pi_{2l-1}^r(\omega):|h(\omega,t)|>0\}\wedge 1, \quad l\in N.
\end{array}
\end{equation*}
\end{remark}

\begin{remark}
\label{rmk:3.5}
Obviously, from the proof of Lemma \ref{lem:3.1}, we can use $\{\tau_k^r\}_{r\in Q^+,k\in N}$ instead of $\{\tau_k^r\}_{k,r\in N}$ in Lemma \ref{lem:3.1}. If we denote $\tau^r:=\lim\limits_{k\rightarrow +\infty}\tau_k^r, \forall r\in Q^+$, for P-a.s. $\omega\in\Omega$, we have
\[
(\bigcup\limits_{r\in {Q^+}}\bigcup\limits_{k\in N}[\tau_k^r(\omega),\tau_{k+1}^r(\omega)])\bigcup(\bigcup\limits_{r\in {Q^+}}\tau^r(\omega))=[0,1], \quad a.e..
\]
In the following, we will use $\{\tau_k^r\}_{r\in Q^+,k\in N}$ instead of $\{\tau_k^r\}_{k,r\in N}$ in Lemma \ref{lem:3.1}.
\end{remark}

Now, we are ready to prove the main result of this paper. Let $(f_n)_{n\geq 0}$ be a sequence of mappings from $[0,1]\times\Omega\times R^{d+d\times m}$ into $R^d$ such that, for all $n\geq 0, f_n=f\ast{\psi}_n$, the convolution product of $f$ and $\psi_n$. By Proposition \ref{pro:2.2}, it follows that $(f_n)_{n\geq 0}$ converges uniformly to $f$ and for any $n\geq 0, f_n$ satisfies Assumption 1 and is Lipschitz with respect to $(y,z)$ uniformly in $(t,\omega)$. In addition, we have
\begin{equation*}
\begin{array}{ll}
|f_n(t,y_1,z)-f_n(t,y_2,z)|\leq \Phi(|y_1-y_2|),\quad \forall y_1,y_2,\\
|f_n(t,y,z_1)-f_n(t,y,z_2)|\leq \Psi(||z_1-z_2||),\quad \forall z_1,z_2.
\end{array}
\end{equation*}
Let $(y_n,z_n)$ be the solution to BSDE associated with $(f_n,\xi)$, that is
\begin{equation}\label{eq:3.2}
\begin{array}{ll}
(y_n,z_n)\in S^{2,d}\times H^{2,d\times m},\\
y_n(t)=\xi+\int_t^1f_n(s,y_n(s),z_n(s))ds-\int_t^1z_n(s)dB_s.
\end{array}
\end{equation}
And from Proposition \ref{pro:2.2}, for any $\varepsilon>0$, there exists $N_{\varepsilon}>0$ such that if $n,m\geq N_{\varepsilon}$, then $|f_n(t,y,z)-f_m(t,y,z)|<\varepsilon$, a.s. In the following, $y_n^i, y^i, \xi^i, f_n^i, f^i$ and $z_n^i, z^i$ denote respectively the $i$th components and rows of $y_n, y, \xi, f_n, f$ and $z_n, z$, $i=1,...,d$, and the linear growth of $\Phi$ and $\Psi$ are denoted by $|\Phi(x)|\leq K(1+|x|)$ and $|\Psi(x)|\leq K(1+|x|), K\geq 0$.\\

\begin{theorem}
\label{thm:3.1}
Suppose that $\xi\in L^2(\Omega,{\cal{F}},P)$ and $f$ satisfies Assumptions 1 and 2, then the following BSDE
\begin{equation}\label{eq:3.3}
y_t=\xi+\int_t^1f(s,y_s,z_s){\rm d}s-\int_t^1z_s{\rm d}B_s,\quad 0\leq t\leq 1,
\end{equation}
has a unique solution $(y,z)$.
\end{theorem}

Without loss of generality, we will suppose $d=2$ in the following.\\

{\bf Proof.}
\textit{Existence.}
\textit{Step 1.} In this step we show that the sequence $(y_n,z_n)_{n\geq 0}$ has a subsequence such that $|y_m-y_n|$ and $E[\int_0^1||z_m(s)-z_n(s)||^2ds]$ are uniformly bounded with respect to $m,n$. When $m,n\geq N_{\varepsilon}$, by It\^{o}'s formula
\begin{equation*}
\begin{array}{lll}
&&\left(y_m \left(t \right)-y_n \left(t \right)\right)^2\\
&=&2\dint_t^1(y_m(s)-y_n(s))(f_m(s,y_m(s),z_m(s))-f_n(s,y_n(s),z_n(s)))ds\\
&&-\dint_t^1(z_m(s)-z_n(s))^2ds-2\dint_t^1(y_m(s)-y_n(s))(z_m(s)-z_n(s))dB_s\\
&=& 2\dint_t^1(y_m(s)-y_n(s))(f_m(s,y_m(s),z_m(s))-f_m(s,y_n(s),z_m(s)))ds\\
&&+2\dint_t^1(y_m(s)-y_n(s))(f_m(s,y_n(s),z_m(s))-f_m(s,y_n(s),z_n(s)))ds\\
&&+2\dint_t^1(y_m(s)-y_n(s))(f_m(s,y_n(s),z_n(s))-f_n(s,y_n(s),z_n(s)))ds\\
&&-\dint_t^1(z_m(s)-z_n(s))^2ds-2\dint_t^1(y_m(s)-y_n(s))(z_m(s)-z_n(s))dB_s\\
&\leq& 2\dint_t^1|y_m(s)-y_n(s)|\Phi(|y_m(s)-y_n(s)|)ds\\
&&+2\dint_t^1|y_m(s)-y_n(s)|\Psi(||z_m(s)-z_n(s)||)ds\\
&&+2\dint_t^1|y_m(s)-y_n(s)||f_m(s,y_n(s),z_n(s))-f_n(s,y_n(s),z_n(s))|ds\\
&&-\dint_t^1(z_m(s)-z_n(s))^2ds-2\dint_t^1(y_m(s)-y_n(s))(z_m(s)-z_n(s))dB_s\\
&\leq& 2K\dint_t^1|y_m(s)-y_n(s)|(|y_m(s)-y_n(s)|+1)ds\\
&&+2K\dint_t^1|y_m(s)-y_n(s)|(||z_m(s)-z_n(s)||+1)ds\\
&&+2\varepsilon\dint_t^1|y_m(s)-y_n(s)|ds-\dint_t^1||z_m(s)-z_n(s)||^2ds\\
&&-2\dint_t^1(y_m(s)-y_n(s))(z_m(s)-z_n(s))dB_s\\
&\leq& ({\varepsilon}^2+4K^2)+(K^2+2K+2)\dint_t^1|y_m(s)-y_n(s)|^2ds\\
&&-2\dint_t^1(y_m(s)-y_n(s))(z_m(s)-z_n(s))dB_s,
\end{array}
\end{equation*}
taking conditional expectation on both sides, for each $0\leq v\leq t\leq 1$, we get
\begin{equation*}
E^{{\cal{F}}_v}[(y_m(t)-y_n(t))^2 ]\leq ({\varepsilon}^2+4K^2)+(K^2+2K+2)\int_t^1E^{{\cal{F}}_v}[|y_m(s)-y_n(s)|^2]ds,
\end{equation*}
by Gronwall's inequality we have $E^{{\cal{F}}_v}[(y_m(t)-y_n(t))^2]\leq ({\varepsilon}^2+4K^2)e^{(K^2+2K+2)}$. If we let $v=t$, we obtain that $(y_m(t)-y_n(t))^2\leq C$, by which we also have $E[\int_0^1||z_m(s)-z_n(s)||^2ds]\leq C$ for $m,n>N_{\varepsilon}$, where $C$ is a constant which may change from one line to another. In the following, for the notational simplicity, the subsequence of $(y_n,z_n)_{n\geq 0}$ is still denoted by $(y_n,z_n)_{n\geq 0}$. \\

\textit{Step 2.} We show that when $(y_n)_{n\geq 0}$ is a Cauchy sequence on $[\tau_1,\tau_2]$ in $S^{2,d}$, it follows also for $(z_n)_{n\geq 0}$ on $[\tau_1,\tau_2]$ in $H^{2,d\times m}$, where $\tau_1$ and $\tau_2$ are two stopping times satisfying $0\leq\tau_1\leq\tau_2\leq 1$ and $P(\tau_1<\tau_2)>0$.\\
Using It\^{o}'s formula, we have on $[\tau_1,\tau_2]$
\begin{equation*}
\begin{array}{lll}
&&E\left[\dint_t^{\tau_2}||z_m(s)-z_n(s)||^2ds\right]\\
&\leq& E\left[(y_m(\tau_2)-y_n(\tau_2))^2\right]-E\left[(y_m(t)-y_n(t))^2\right]\\
&&+2E\left[\dint_t^{\tau_2}(y_m(s)-y_n(s))(f_m(s,y_m(s),z_m(s))-f_n(s,y_n(s),z_n(s)))ds\right]\\
&\leq& E\left[(y_m(\tau_2)-y_n(\tau_2))^2\right]-E\left[(y_m(t)-y_n(t))^2\right]\\
&&+2\sqrt{E\left[\sup\limits_{\tau_1\leq s\leq \tau_2}|y_m(s)-y_n(s)|^2\right]}\\
&&\times\sqrt{E\left[\dint_t^{\tau_2}|f_m(s,y_m(s),z_m(s))-f_n(s,y_n(s),z_n(s))|^2ds\right]},
\end{array}
\end{equation*}
On the other hand, Step 1 implies that there exists a constant $C\geq 0$ such that, for all $m,n\geq 0$,
\[
E\left[\int_t^{\tau_2}|f_m(s,y_m(s),z_m(s))-f_n(s,y_n(s),z_n(s))|^2ds\right]\leq C.
\]
Hence, the sequence $(z_n)_{n\geq 0}$ is a Cauchy sequence on $[\tau_1,\tau_2]$ in $H^{2,d\times m}$.\\

\textit{Step 3.} \textit{(i)} For all $m,n, m\neq n$, we denote $z_{mn}(t):=\sum\limits_{i=1}^2\text{sgn}(y_m^i(t)-y_n^i(t))(z_m^i(t)-z_n^i(t))$ and $Z_{mn}(t):=\int_0^t|z_{mn}(s)|ds$, respectively. By Lemma \ref{lem:3.1}, we know that there exist countable stopping times $\{\tau_k^{r,mn}\}_{r\in Q^+, k\in N}$ such that, for P-a.s. $\omega\in\{\omega:\tau_k^{r,mn}(\omega)<\tau_{k+1}^{r,mn}(\omega)\}$, $Z_{mn}(\omega,\cdot)$ is strictly increasing on $[\tau_k^{r,mn}(\omega),\tau_{k+1}^{r,mn}(\omega)]$ or $Z_{mn}(\omega,\tau_k^{r,mn}(\omega))=Z_{mn}(\omega,\tau_{k+1}^{r,mn}(\omega))$. Moreover, for P-a.s. $\omega\in \Omega$,  $[0,1]=\\
\bigcup\limits_{r\in Q^+}\bigcup\limits_{k\in N} [\tau_k^{r,mn}(\omega),\tau_{k+1}^{r,mn}(\omega)]
$, a.e. In the following, as Remark \ref{rmk:3.4}, for given $m,n, m\neq n$, we will take a version of the process $z_{mn}$ such that, for P-a.s. $\omega\in\{\omega:\tau_k^{r,mn}(\omega)<\tau_{k+1}^{r,mn}(\omega)\}$, $|z_{mn}(\omega,\cdot)|>0$ on $(\tau_k^{r,mn}(\omega),\tau_{k+1}^{r,mn}(\omega))$ or $|z_{mn}(\omega,\cdot)|=0$ on $[\tau_k^{r,mn}(\omega),\tau_{k+1}^{r,mn}(\omega)]$.\\

\textit{(ii)} For given $r_0\in Q^+, k_0\geq 1$ satisfying $P(\tau_{k_0}^{r_0, mn}<\tau_{k_0+1}^{r_0, mn})>0$ and $\forall m,n>N_{\varepsilon}, m\neq n$, if $|z_{mn}|>0$ on $[\tau_{k_0}^{r_0, mn},\tau_{k_0+1}^{r_0, mn}]$, for all $\varepsilon_0\in (0,1)$, we denote
\[
{\tau}_{\varepsilon_0}^{mn}(\omega):=\inf\left\{t>\tau_{k_0}^{r_0,mn}(\omega):|z_{mn}(t,\omega)|\geq\varepsilon_0\right\}\wedge \tau_{k_0+1}^{r_0,mn},
\]
\[
{\vartheta}_{\varepsilon_0}^{mn}(\omega):=\inf\left\{t>\tau_{k_0}^{r_0,mn}(\omega):||z_m(t,\omega)-z_n(t,\omega)||\geq\frac{1}{\varepsilon_0}\right\}\wedge \tau_{k_0+1}^{r_0,mn},
\]
\begin{equation*}
{\nu}_{\varepsilon_0}^{mn}(\omega):=\inf\left\{t>{\tau}_{\varepsilon_0}^{mn}(\omega):|z_{mn}(s,\omega)|<\varepsilon_0\right\}\wedge \tau_{k_0+1}^{r_0,mn},
\end{equation*}
\[
{\tau}_{\varepsilon_0,1}^{mn}:={\nu}_{\varepsilon_0}^{mn}\wedge{\vartheta}_{\varepsilon_0}^{mn},\quad
{\tau}_{\varepsilon_0,2}^{mn}:={\tau}_{\varepsilon_0}^{mn}\wedge {\tau}_{\varepsilon_0,1}^{mn}.
\]
Obviously, ${\tau}_{\varepsilon_0}^{mn}, {\vartheta}_{\varepsilon_0}^{mn}$ and ${\nu}_{\varepsilon_0}^{mn}$ are all stopping times, this implies that ${\tau}_{\varepsilon_0,1}^{mn}$ and ${\tau}_{\varepsilon_0,2}^{mn}$ are also stopping times. Since $P(\tau_{k_0}^{r_0, mn}<\tau_{k_0+1}^{r_0, mn})>0$ and P-a.s. $|z_{mn}|>0$ on $\{\omega:\tau_{k_0}^{r_0,mn}(\omega)<\tau_{k_0+1}^{r_0,mn}(\omega)\}$, then we have P-a.s. ${\tau}_{\varepsilon_0,2}^{mn}\rightarrow \tau_{k_0}^{r_0,mn}$, ${\tau}_{\varepsilon_0,1}^{mn}\rightarrow \tau_{k_0+1}^{r_0,mn}$ as $\varepsilon_0\rightarrow 0$. Hence, there exists $\varepsilon_0\in (0,1)$ such that $P({\tau}_{\varepsilon_0,2}^{mn}<{\tau}_{\varepsilon_0,1}^{mn})>0$.\\

\textit{(iii)} For $i\in\{1,2\}$, by Tanaka's formula, we have
\begin{equation*}
\begin{array}{lll}
&&|y_n^i(t)-y_m^i(t)|\\
&\leq& |y_n^i({\tau}_{\varepsilon_0,1}^{mn})-y_m^i({\tau}_{\varepsilon_0,1}^{mn})|\\
&&+\dint_t^{{\tau}_{\varepsilon_0,1}^{mn}}
\text{sgn}(y_n^i(s)-y_m^i(s))(f_n^i(s,y_n(s),z_n(s))-f_m^i(s,y_m(s),z_m(s)))ds\\
&&-\dint_t^{{\tau}_{\varepsilon_0,1}^{mn}}\text{sgn}(y_n^i(s)-y_m^i(s))(z_n^i(s)-z_m^i(s))dB_s\\
&\leq&|y_n^i({\tau}_{\varepsilon_0,1}^{mn})-y_m^i({\tau}_{\varepsilon_0,1}^{mn})|+\dint_t^{{\tau}_{\varepsilon_0,1}^{mn}}(\Phi(|y_n(s)-y_m(s)|)+\varepsilon)ds\\
&&+\dint_t^{{\tau}_{\varepsilon_0,1}^{mn}}\Psi(||z_n(s)-z_m(s)||)ds\\
&&-\dint_t^{{\tau}_{\varepsilon_0,1}^{mn}}\text{sgn}(y_n^i(s)-y_m^i(s))(z_n^i(s)-z_m^i(s))dB_s, \quad {\tau}_{\varepsilon_0,2}^{mn}\leq t\leq {\tau}_{\varepsilon_0,1}^{mn},
\end{array}
\end{equation*}
it follows that
\begin{equation*}
\begin{array}{lll}
&&\sum\limits_{i=1}^2|y_n^i(t)-y_m^i(t)|\\
&\leq&\sum\limits_{i=1}^2|y_n^i({\tau}_{\varepsilon_0,1}^{mn})-y_m^i({\tau}_{\varepsilon_0,1}^{mn})|+3\dint_t^{{\tau}_{\varepsilon_0,1}^{mn}}(\Phi(|y_n(s)-y_m(s)|)+\varepsilon)ds\\
&&+2\dint_t^{{\tau}_{\varepsilon_0,1}^{mn}}\Psi(||z_n(s)-z_m(s)||)ds\\
&&-\dint_t^{{\tau}_{\varepsilon_0,1}^{mn}}\sum\limits_{i=1}^2\text{sgn}(y_n^i(s)-y_m^i(s))(z_n^i(s)-z_m^i(s))dB_s\\
&\leq&\sum\limits_{i=1}^2|y_n^i({\tau}_{\varepsilon_0,1}^{mn})-y_m^i({\tau}_{\varepsilon_0,1}^{mn})|+3\dint_t^{{\tau}_{\varepsilon_0,1}^{mn}}(\Phi(|y_n(s)-y_m(s)|)+\varepsilon)ds\\
&&-\dint_t^{{\tau}_{\varepsilon_0,1}^{mn}}\sum\limits_{i=1}^2\text{sgn}(y_n^i(s)-y_m^i(s))(z_n^i(s)-z_m^i(s))dB_s^{mn}, \quad {\tau}_{\varepsilon_0,2}^{mn}\leq t\leq {\tau}_{\varepsilon_0,1}^{mn},
\end{array}
\end{equation*}
the process $B^{mn}$ is a Brownian motion on $[0,{\tau}_{\varepsilon_0,1}^{mn}]$ under $P^{mn}$ which is equivalent to $P$ and defined by
\[
\frac{dP^{mn}}{dP}=\exp\left(\int_0^{{\tau}_{\varepsilon_0,1}^{mn}}\eta_{mn}(s)dB_s
-\frac{1}{2}\int_0^{{\tau}_{\varepsilon_0,1}^{mn}}|\eta_{mn}(s)|^2ds\right),
\]
where
\begin{equation*}
\eta_{mn}(t):=
\left\{
\begin{array}{lll}
\frac{2\Psi(||z_n(t)-z_m(t)||)}{\sum\limits_{i=1}^2\text{sgn}(y_n^i(t)-y_m^i(t))(z_n^i(t)-z_m^i(t))},\quad {\tau}_{\varepsilon_0,2}^{mn}\leq t\leq {\tau}_{\varepsilon_0,1}^{mn},\\
0,\qquad\qquad 0\leq t<{\tau}_{\varepsilon_0,2}^{mn}.
\end{array}
\right.
\end{equation*}
Then taking the conditional expectation under $P^{mn}$, we have
\begin{equation}\label{eq:3.4}
\begin{array}{lll}
&&\sum\limits_{i=1}^2|y_m^i(t)-y_n^i(t)|\\
&\leq&E^{mn}\left[\sum\limits_{i=1}^2|y_m^i({\tau}^{mn}_{\varepsilon_0,1})-y_n^i({\tau}^{mn}_{\varepsilon_0,1})||{\cal{F}}_t\right]\\
&&+E^{mn}\left[3\dint_t^{{\tau}_{\varepsilon_0,1}}(\Phi(|y_m(s)-y_n(s)|)+\varepsilon)ds|{\cal{F}}_t\right],\quad {\tau}_{\varepsilon_0,2}^{mn}\leq t\leq {\tau}_{\varepsilon_0,1}^{mn}.
\end{array}
\end{equation}

Denoting $\Phi_{\varepsilon}(x):=3(\Phi(x)+\varepsilon), x\in R^+$. By the continuity of $y_m, y_n$ and $|y_m(t)-y_n(t)|\leq C, \forall m,n, m\neq n$, we know that, for all stopping time $\vartheta$ such that $0\leq \vartheta\leq 1$, a.s., there exists a $\Phi_{\varepsilon}$-envelope of the process $|y_m^1-y_n^1|+|y_m^2-y_n^2|$ at $\vartheta$ denoted by $\theta(\vartheta)$ such that $|y_m^1(\vartheta)-y_n^1(\vartheta)|+|y_m^2(\vartheta)-y_n^2(\vartheta)|\leq \theta(\vartheta)$, a.s.. Moreover, $\theta(\vartheta)\rightarrow 0$ as $\vartheta\rightarrow 1$.\\
Now, we consider the following equation:
\begin{equation}\label{eq:3.5}
u^{\varepsilon_0}(t)=\theta(\tau_{\varepsilon_0,1}^{mn})
+3\dint_t^{\tau_{\varepsilon_0,1}^{mn}}(\Phi(u^{\varepsilon_0}(s))+\varepsilon)ds,
\end{equation}
where $0\leq t\leq \tau_{\varepsilon_0,1}^{mn}$. From the definition of $\theta(\tau_{\varepsilon_0,1}^{mn})$, we know Eq. (\ref{eq:3.5}) has a unique solution $u^{\varepsilon_0}(t)$ and $u^{\varepsilon_0}(t), \forall t,$ is a constant. Taking the conditional expectation on both side, it follows that
\[
u^{\varepsilon_0}(t)=E^{mn}[\theta(\tau_{\varepsilon_0,1}^{mn})
+3\dint_t^{\tau_{\varepsilon_0,1}^{mn}}(\Phi(u^{\varepsilon_0}(s))+\varepsilon)ds|{\cal{F}}_t].
\]
From Appendix, we have P-a.s.
\begin{equation}\label{eq:3.6}
|y_m^1(t)-y_n^1(t)|+|y_m^2(t)-y_n^2(t)|\leq u^{\varepsilon_0}(t),\quad {\tau}_{\varepsilon_0,2}^{mn}\leq t\leq {\tau}_{\varepsilon_0,1}^{mn}.
\end{equation}
As $\varepsilon_0\rightarrow 0$, we have ${\tau}_{\varepsilon_0,2}^{mn}\rightarrow \tau^{r_0,mn}_{k_0}, {\tau}_{\varepsilon_0,1}^{mn}\rightarrow \tau^{r_0,mn}_{k_0+1}$ and $|y_m^i({\tau}_{\varepsilon_0,1}^{mn})-y_n^i({\tau}_{\varepsilon_0,1}^{mn})|\rightarrow |y_m^i(\tau^{r_0,mn}_{k_0+1})-y_n^i(\tau^{r_0,mn}_{k_0+1})|, i=1,2$, this implies that P-a.s.
\[
\lim\limits_{\varepsilon_0\rightarrow 0}\theta(\tau_{\varepsilon_0,1}^{mn})=\theta(\tau^{r_0,mn}_{k_0+1}),
\]
and
\begin{equation*}
\lim\limits_{\varepsilon_0\rightarrow 0}u^{\varepsilon_0}(t)= u^{r_0}_{k_0+1}(t),
\end{equation*}
where
\begin{equation*}
u^{r_0}_{k_0+1}(t)=\theta(\tau^{r_0,mn}_{k_0+1})
+3\dint_t^{\tau^{r_0,mn}_{k_0+1}}(\Phi(u^{r_0}_{k_0+1}(s))+\varepsilon)ds.
\end{equation*}
Hence, as $\varepsilon_0\rightarrow 0$, we can obtain P-a.s.
\begin{equation}\label{eq:3.7}
|y_m^1(t)-y_n^1(t)|+|y_m^2(t)-y_n^2(t)|\leq u^{r_0}_{k_0+1}(t), \quad\tau^{r_0,mn}_{k_0}\leq t\leq \tau^{r_0,mn}_{k_0+1}.
\end{equation}

\textit{Step 4}. \textit{(i)} If $|z_{mn}|=0$ on $[\tau^{r_0,mn}_{k_0},\tau^{r_0,mn}_{k_0+1}]$, we denote $Z_{mn}^1(t):=\int_{\tau^{r_0,mn}_{k_0}}^t|z_m^1(s)-z_n^1(s)|ds, \tau^{r_0,mn}_{k_0}\leq t\leq \tau^{r_0,mn}_{k_0+1}$.  By Lemma \ref{lem:3.1}, it follows that there exist countable stopping times $\{\rho_k^{r,mn}\}_{r\in Q^+,k\in N}$ such that, for P-a.s.  $\omega\in\{\omega:\rho_k^{r,mn}(\omega)<\rho_{k+1}^{r,mn}(\omega)\}$, $Z_{mn}^1(\omega,\rho_k^{r,mn}(\omega))=Z_{mn}^1(\omega,\rho_{k+1}^{r,mn}(\omega))$ or $Z_{mn}^1(\omega,\cdot)$ is strictly increasing on $[\rho_k^{r,mn}(\omega),\rho_{k+1}^{r,mn}(\omega)]$. Moreover, for P-a.s. $\omega\in \Omega$,  $\bigcup\limits_{r\in Q^+ }\bigcup\limits_{k\in N }[\rho_k^{r,mn}(\omega),\rho_{k+1}^{r,mn}(\omega)]=[\tau^{r_0,mn}_{k_0}(\omega),\tau^{r_0,mn}_{k_0+1}(\omega)]$, a.e.. In the following, as Remark \ref{rmk:3.4}, for given $r$ and $k$, we will take a version of $z_m$ such that $z_m^1(t)-z_n^1(t)=z_m^2(t)-z_n^2(t)=0$ on $[\rho_k^{r,mn},\rho_{k+1}^{r,mn}]$ or $\text{sgn}(y_m^1(t)-y_n^1(t))(z_m^1(t)-z_n^1(t))=-\text{sgn}(y_m^2(t)-y_n^2(t))(z_m^2(t)-z_n^2(t))\neq 0$ on $(\rho_k^{r,mn},\rho_{k+1}^{r,mn})$.\\

\textit{(ii)} For given $r_1$ and $k_1$ satisfying $P(\rho_{k_1}^{r_1,mn}<\rho_{k_1+1}^{r_1,mn})>0$, we first consider the case of $z_m^1-z_n^1=z_m^2-z_n^2=0$ on $[\rho_{k_1}^{r_1,mn},\rho_{k_1+1}^{r_1,mn}]$. By Tanaka's formula, we have
\[
\sum\limits_{i=1}^2|y_n^i(t)-y_m^i(t)|\leq\sum\limits_{i=1}^2|y_n^i(\rho_{k_1+1}^{r_1,mn})-y_m^i(\rho_{k_1+1}^{r_1,mn})|
+3\dint_t^{\rho_{k_1+1}^{r_1,mn}}(\Phi(|y_n(s)-y_m(s)|)+\varepsilon)ds.\\
\]
Obviously, we have
\begin{equation}\label{eq:3.8}
\sum\limits_{i=1}^2|y_n^i(t)-y_m^i(t)|\leq u_{k_1+1}^{r_1}(t), \quad \rho_{k_1}^{r_1,mn}\leq t\leq\rho_{k_1+1}^{r_1,mn},
\end{equation}
where $u_{k_1+1}^{r_1}(t)$ is the solution of the following equation:
\[
u^{r_1}_{k_1+1}(t)=\theta(\rho^{r_1,mn}_{k_1+1})
+3\dint_t^{\rho^{r_1,mn}_{k_1+1}}(\Phi(u^{r_1}_{k_1+1}(s))+\varepsilon)ds,
\]
and $\theta(\rho^{r_1,mn}_{k_1+1})$ is the $\Phi_{\varepsilon}$-envelope of the process $\sum\limits_{i=1}^2|y_n^i-y_m^i|$ at $\rho^{r_1,mn}_{k_1+1}$.\\

\textit{(iii)} If $\text{sgn}(y_m^1(t)-y_n^1(t))(z_m^1(t)-z_n^1(t))=-\text{sgn}(y_m^2(t)-y_n^2(t))(z_m^2(t)-z_n^2(t))\neq 0$ on $(\rho_{k_1}^{r_1,mn},\rho_{k_1+1}^{r_1,mn})$. For all $\varepsilon_0\in (0,1)$, we denote
\[
{\tau}_{\varepsilon_0}^{mn,1}(\omega):=\inf\left\{t>\rho^{r_1,mn}_{k_1}(\omega):|z_m^1(t,\omega)-z_n^1(t,\omega)|\geq\varepsilon_0\right\}\wedge \rho^{r_1,mn}_{k_1+1},
\]
\[
{\vartheta}_{\varepsilon_0,1}^{mn}(\omega):=\inf\left\{t>\rho^{r_1,mn}_{k_1}(\omega):||z_m(t,\omega)-z_n(t,\omega)||\geq\frac{1}{\varepsilon_0}\right\}\wedge \rho^{r_1,mn}_{k_1+1},
\]
\begin{equation*}
{\nu}_{\varepsilon_0}^{mn,1}(\omega):=\inf\left\{t>{\tau}_{\varepsilon_0}^{mn,1}(\omega):|z_m^1(t,\omega)-z_n^1(t,\omega)|<\varepsilon_0\right\}\wedge \rho^{r_1,mn}_{k_1+1},
\end{equation*}
\[
{\tau}_{\varepsilon_0,1}^{mn,1}:={\nu}_{\varepsilon_0}^{mn,1}\wedge{\vartheta}_{\varepsilon_0,1}^{mn},\quad
{\tau}_{\varepsilon_0,2}^{mn,1}:={\tau}_{\varepsilon_0}^{mn,1}\wedge {\tau}_{\varepsilon_0,1}^{mn,1}.
\]
Obviously, we have P-a.s. ${\tau}_{\varepsilon_0,2}^{mn,1}\rightarrow \rho_{k_1}^{r_1,mn}$, ${\tau}_{\varepsilon_0,1}^{mn,1}\rightarrow \rho_{k_1+1}^{r_1,mn}$ as $\varepsilon_0\rightarrow 0$. Hence, there exists $\varepsilon_0\in (0,1)$ such that $P({\tau}_{\varepsilon_0,2}^{mn,1}<{\tau}_{\varepsilon_0,1}^{mn,1})>0$.\\
Then by Tanaka's formula, it follows that, $\forall \alpha\in(0,1)$,
\begin{equation*}
\begin{array}{lll}
&&(1+\alpha)|y_n^1(t)-y_m^1(t)|+|y_n^2(t)-y_m^2(t)|\\
&\leq&(1+\alpha)|y_n^1({\tau}_{\varepsilon_0,1}^{mn,1})-y_m^1({\tau}_{\varepsilon_0,1}^{mn,1})|
+|y_n^2({\tau}_{\varepsilon_0,1}^{mn,1})-y_m^2({\tau}_{\varepsilon_0,1}^{mn,1})|\\
&&+3\dint_t^{{\tau}_{\varepsilon_0,1}^{mn,1}}(\Phi(|y_n(s)-y_m(s)|)+\varepsilon)ds
+3\dint_t^{{\tau}_{\varepsilon_0,1}^{mn,1}}\Psi(||z_n(s)-z_m(s)||)ds\\
&&-\alpha\dint_t^{{\tau}_{\varepsilon_0,1}^{mn,1}}\text{sgn}(y_n^1(s)-y_m^1(s))(z_n^1(s)-z_m^1(s))dB_s\\
&=&(1+\alpha)|y_n^1({\tau}_{\varepsilon_0,1}^{mn,1})-y_m^1({\tau}_{\varepsilon_0,1}^{mn,1})|
+|y_n^2({\tau}_{\varepsilon_0,1}^{mn,1})-y_m^2({\tau}_{\varepsilon_0,1}^{mn,1})|\\
&&+3\dint_t^{{\tau}_{\varepsilon_0,1}^{mn,1}}(\Phi(|y_n(s)-y_m(s)|)+\varepsilon)ds\\
&&-\alpha\dint_t^{{\tau}_{\varepsilon_0,1}^{mn,1}}\text{sgn}(y_n^1(s)-y_m^1(s))(z_n^1(s)-z_m^1(s))dB_s^{mn,1},
\quad {\tau}_{\varepsilon_0,2}^{mn,1}\leq t\leq{\tau}_{\varepsilon_0,1}^{mn,1},\\
\end{array}
\end{equation*}
the process $B^{mn,1}$ is a Brownian motion on $[0,{\tau}_{\varepsilon_0,1}^{mn,1}]$ under $P^{mn,1}$ which is equivalent to $P$ and defined by
\[
\frac{dP^{mn,1}}{dP}=\exp\left(\int_0^{{\tau}_{\varepsilon_0,1}^{mn,1}}\eta_{mn,1}(s)dB_s
-\frac{1}{2}\int_0^{{\tau}_{\varepsilon_0,1}^{mn,1}}|\eta_{mn,1}(s)|^2ds\right),
\]
where
\begin{equation*}
\eta_{mn,1}(t):=
\left\{
\begin{array}{lll}
\frac{3\Psi(||z_n(t)-z_m(t)||)}{\alpha\cdot\text{sgn}(y_n^1(t)-y_m^1(t))(z_n^1(t)-z_m^1(t))},\quad {\tau}_{\varepsilon_0,2}^{mn,1}\leq t\leq {\tau}_{\varepsilon_0,1}^{mn,1},\\
0,\qquad\qquad 0\leq t<{\tau}_{\varepsilon_0,2}^{mn,1}.
\end{array}
\right.
\end{equation*}
Then taking the conditional expectation under $P^{mn,1}$, we have
\begin{equation}\label{eq:3.9}
\begin{array}{lll}
&&(1+\alpha)|y_n^1(t)-y_m^1(t)|+|y_n^2(t)-y_m^2(t)|\\
&\leq&E^{mn,1}\left[(1+\alpha)|y_n^1({\tau}_{\varepsilon_0,1}^{mn,1})-y_m^1({\tau}_{\varepsilon_0,1}^{mn,1})|
+|y_n^2({\tau}_{\varepsilon_0,1}^{mn,1})-y_m^2({\tau}_{\varepsilon_0,1}^{mn,1})||{\cal{F}}_t\right]\\
&&+E^{mn,1}\left[3\dint_t^{{\tau}_{\varepsilon_0,1}^{mn,1}}(\Phi(|y_m(s)-y_n(s)|)+\varepsilon)ds|{\cal{F}}_t\right],\quad {\tau}_{\varepsilon_0,2}^{mn,1}\leq t\leq {\tau}_{\varepsilon_0,1}^{mn,1}.
\end{array}
\end{equation}
Let $\theta_{\alpha}({\tau}_{\varepsilon_0,1}^{mn,1})$ be the $\Phi_{\varepsilon}$-envelope of the process $(1+\alpha)|y_n^1(t)-y_m^1(t)|+|y_n^2(t)-y_m^2(t)|$ at ${\tau}_{\varepsilon_0,1}^{mn,1}$ and we consider the following equation:
\[
u^{\varepsilon_0,1}_{\alpha}(t)=\theta_{\alpha}({\tau}_{\varepsilon_0,1}^{mn,1})
+3\dint_t^{{\tau}_{\varepsilon_0,1}^{mn,1}}(\Phi(u^{\varepsilon_0,1}_{\alpha}(s))+\varepsilon)ds.
\]
Obviously, we have
\[
(1+\alpha)|y_n^1(t)-y_m^1(t)|+|y_n^2(t)-y_m^2(t)|\leq u^{\varepsilon_0,1}_{\alpha}(t),\quad  {\tau}_{\varepsilon_0,2}^{mn,1}\leq t\leq {\tau}_{\varepsilon_0,1}^{mn,1}.
\]
First, let $\varepsilon_0\rightarrow 0$, it follows that
\begin{equation}\label{eq:3.10}
(1+\alpha)|y_n^1(t)-y_m^1(t)|+|y_n^2(t)-y_m^2(t)|\leq u^{r_1}_{k_1+1,\alpha}(t), \quad\rho^{r_1,mn}_{k_1}\leq t\leq \rho^{r_1,mn}_{k_1+1},
\end{equation}
where
\[
u^{r_1}_{k_1+1,\alpha}(t)=\theta_{\alpha}(\rho^{r_1,mn}_{k_1+1})
+3\dint_t^{\rho^{r_1,mn}_{k_1+1}}(\Phi(u^{r_1}_{k_1+1,\alpha}(s))+\varepsilon)ds.
\]
On the other hand, from the fact that $\sum\limits_{i=1}^2|y_n^i(t)-y_m^i(t)|\leq (1+\alpha)|y_n^1(t)-y_m^1(t)|+|y_n^2(t)-y_m^2(t)|\leq (1+\alpha)\sum\limits_{i=1}^2|y_n^i(t)-y_m^i(t)|$, we can conclude that $\theta_{\alpha}(\rho^{r_1,mn}_{k_1+1})\rightarrow \theta(\rho^{r_1,mn}_{k_1+1})$ as $\alpha\rightarrow 0$. Hence, as $\alpha\rightarrow 0$, it follows that $u^{r_1}_{k_1+1,\alpha}(t)\rightarrow u^{r_1}_{k_1+1}(t)$ and
\begin{equation}\label{eq:3.11}
\sum\limits_{i=1}^2|y_n^i(t)-y_m^i(t)|\leq u_{k_1+1}^{r_1}(t), \quad \rho_{k_1}^{r_1,mn}\leq t\leq\rho_{k_1+1}^{r_1,mn},
\end{equation}
where
\[
u^{r_1}_{k_1+1}(t)=\theta(\rho^{r_1,mn}_{k_1+1})
+3\dint_t^{\rho^{r_1,mn}_{k_1+1}}(\Phi(u^{r_1}_{k_1+1}(s))+\varepsilon)ds.
\]\\

\textit{(iv)} Let $u_{l+1}^{r_1}(t)$ be the solution to the following equation:
\begin{equation}\label{eq:3.12}
u^{r_1}_{l+1}(t)=\theta(\rho^{r_1,mn}_{l+1})
+3\dint_t^{\rho^{r_1,mn}_{l+1}}(\Phi(u^{r_1}_{l+1}(s))+\varepsilon)ds,
\end{equation}
where $l>k_1+1$, then we also have
\begin{equation}\label{eq:3.13}
\sum\limits_{i=1}^2|y_n^i(t)-y_m^i(t)|\leq u_{l+1}^{r_1}(t), \quad \rho_{k_1}^{r_1,mn}\leq t\leq\rho_{k_1+1}^{r_1,mn}.
\end{equation}
In fact, similar to the proof of the inequality (\ref{eq:3.8}) and (\ref{eq:3.11}), we have P-a.s.
\begin{equation}\label{eq:3.14}
\sum\limits_{i=1}^2|y_n^i(t)-y_m^i(t)|\leq u_{l+1}^{r_1}(t), \quad \rho_{l}^{r_1,mn}\leq t\leq\rho_{l+1}^{r_1,mn}.
\end{equation}
Then, on $[0,\rho_{l}^{r_1,mn}]$, Eq. (\ref{eq:3.12}) becomes
\[
u^{r_1}_{l+1}(t)=u^{r_1}_{l+1}(\rho^{r_1,mn}_l)
+3\dint_t^{\rho^{r_1,mn}_l}(\Phi(u^{r_1}_{l+1}(s))+\varepsilon)ds,
\]
By (\ref{eq:3.14}) and the definition of $\theta(\rho^{r_1,mn}_l)$,  we have $\theta(\rho^{r_1,mn}_l)\leq u^{r_1}_{l+1}(\rho^{r_1,mn}_l)$. Hence, on $[\rho_{l-1}^{r_1,mn},\rho_l^{r_1,mn}]$, we have
\[
\sum\limits_{i=1}^2|y_n^i(t)-y_m^i(t)|\leq u_{l+1}^{r_1}(t).
\]
It follows that
\[
\sum\limits_{i=1}^2|y_n^i(t)-y_m^i(t)|\leq u_{l+1}^{r_1}(t), \quad \rho_{l-1}^{r_1,mn}\leq t\leq\rho_{l+1}^{r_1,mn}.
\]
Taking the above procedure step by step, we can obtain
\[
\sum\limits_{i=1}^2|y_n^i(t)-y_m^i(t)|\leq u_{l+1}^{r_1}(t), \quad \rho_{k_1}^{r_1,mn}\leq t\leq\rho_{l+1}^{r_1,mn},
\]
and we get the inequality (\ref{eq:3.13}). \\
If we let $l\rightarrow +\infty$ and denote $\rho^{r_1,mn}:=\lim\limits_{l\rightarrow +\infty}\rho_l^{r_1,mn}$, it is easy to obtain that
\[
\sum\limits_{i=1}^2|y_n^i(t)-y_m^i(t)|\leq u^{r_1}(t), \quad \rho_{k_1}^{r_1,mn}\leq t\leq\rho^{r_1,mn},
\]
where
\[
u^{r_1}(t)=\theta(\rho^{r_1,mn})
+3\dint_t^{\rho^{r_1,mn}}(\Phi(u^{r_1}(s))+\varepsilon)ds.
\]

\textit{(v)} For all stopping time sequence $\{\rho_k^{r,mn}\}_{k\in N}$, we denote $\rho^{r,mn}:=\lim\limits_{k\rightarrow +\infty}\rho_k^{r,mn}$. Now, we show that, for all $r_2>r_1$ and  $\forall\rho_{k_2+1}^{r_2,mn}\in \{\rho^{r_2,mn}_k\}_{k\in N}$, it still holds that
\begin{equation}\label{eq:3.15}
\sum_{i=1}^2|y_m^i(t,\omega)-y_n^i(t,\omega)|\leq u^{r_2}_{k_2+1}(t,\omega), \quad\rho^{r_1,mn}_{k_1}(\omega)\leq t\leq \rho^{r_1,mn}_{k_1+1}(\omega),
\end{equation}
for P-a.s. $\omega\in\{\omega:\rho_{k_2+1}^{r_2,mn}(\omega)\geq\rho^{r_1,mn}_{k_1+1}(\omega)\}$, where
\begin{equation*}
u^{r_2}_{k_2+1}(t)=\theta(\rho^{r_2,mn}_{k_2+1})+3\dint_t^{\rho^{r_2,mn}_{k_2+1}}(\Phi(u^{r_2}_{k_2+1}(s))+\varepsilon)ds.
\end{equation*}
Similar to (ii) and (iii), we can obtain
\begin{equation}\label{eq:3.16}
\sum_{i=1}^2|y_m^i(t)-y_n^i(t)|\leq u^{r_2}_{k_2+1}(t), \quad \rho^{r_2,mn}_{k_2}\leq t\leq \rho^{r_2,mn}_{k_2+1},\quad a.s.
\end{equation}
Let ${\cal{A}}=\{\rho^{r,mn}_k: r\in Q_1, k\in N\}$ be a set of all stopping times, where $Q_1\subseteq[r_1,r_2]\cap Q^+$, such that
\[
\sum_{i=1}^2|y_m^i(t)-y_n^i(t)|\leq u^{r_2}_{k_2+1}(t), \quad \rho^{r,mn}_k\leq t\leq \rho^{r_2,mn}_{k_2+1},\quad a.s..
\]
Obviously, $\rho^{r_2,mn}_{k_2}\in {\cal{A}}$. Hence, ${\cal{A}}$ is not an empty set. On the other hand, similar to (iv), we know if there exist $r_3\in [r_1,r_2)\cap Q^+$ and $k_3\in N$ such that
\[
\sum_{i=1}^2|y_m^i(t)-y_n^i(t)|\leq u^{r_2}_{k_2+1}(t), \quad \rho^{r_3,mn}_{k_3}\leq t\leq \rho^{r_2,mn}_{k_2+1},\quad a.s.,
\]
then for all $k\in N$, it still holds that
\[
\sum_{i=1}^2|y_m^i(t)-y_n^i(t)|\leq u^{r_2}_{k_2+1}(t), \quad \rho^{r_3,mn}_k\leq t\leq \rho^{r_2,mn}_{k_2+1},\quad a.s..
\]
Now, we denote $\rho:=\text{essinf}\{\rho^{r,mn}_k: r\in Q_1, k\in N\}$. For all $\rho^{r_i,mn}_{k_i}, \rho^{r_j,mn}_{k_j}\in {\cal{A}}$, it obviously holds that $\rho^{r_i,mn}_{k_i}\wedge\rho^{r_j,mn}_{k_j}\in {\cal{A}}$. And by Lemma \ref{lem:2.1}, we know $\rho$ is a stopping time. By the continuity of $|y_m(\cdot)-y_n(\cdot)|$ and $u^{r_1}_{k_1+1}(\cdot)$, it follows that
\[
\sum_{i=1}^2|y_m^i(t)-y_n^i(t)|\leq u^{r_2}_{k_2+1}(t), \quad \rho\leq t\leq \rho^{r_2,mn}_{k_2+1},\quad a.s.
\]
Now, we show $\rho\leq\rho^{r_1,mn}_{k_1}$, a.s.. Otherwise, we denote $Q_2:=([r_1,r_2]\cap Q^+)\setminus Q_1$. By the continuity of $Z_{mn}^1$ and the fact that, for P-a.s. $\omega\in \Omega$,
\[
(\bigcup\limits_{r\in {Q^+}}\bigcup\limits_{k\in N}[\rho_k^{r,mn}(\omega),\rho_{k+1}^{r,mn}(\omega)])\bigcup(\bigcup\limits_{r\in {Q^+}}\rho^{r,mn}(\omega))=[\tau^{r_0,mn}_{k_0}(\omega),\tau^{r_0,mn}_{k_0+1}(\omega)],\quad a.e.,
\]
it is easy to obtain that there exists $r_4\in Q_2$ such that $\rho^{r_4,mn}\geq\rho$, a.s.. By (iv), it follows that P-a.s.
\[
\sum_{i=1}^2|y_m^i(t)-y_n^i(t)|\leq u^{r_2}_{k_2+1}(t), \quad \rho^{r_4,mn}_k\leq t\leq \rho^{r_2,mn}_{k_2+1},\quad \forall k\in N,
\]
this contradicts to the definition of $Q_1$. Therefore it holds that $\rho\leq\rho^{r_1,mn}_{k_1}$ and we get the inequality (\ref{eq:3.15}).\\

\textit{(vi)} Now, we take a subsequence of $\{\rho^{r,mn}_k\}_{r\in Q^+,k\in N}$ denoted by $\{\rho^{mn}_k\}_{k\in N}$ such that $\rho^{mn}_k\rightarrow \tau^{r_0,mn}_{k_0+1}$ as $k\rightarrow +\infty$ and we consider the following equation:
\begin{equation*}
v_k(t)=\theta(\rho^{mn}_k)+3\int_t^{\rho^{mn}_k}(\Phi(v_k(s))+\varepsilon)ds,\\
\end{equation*}
then we have P-a.s. $\sum_{i=1}^2|y_m^i(t)-y_n^i(t)|\leq v_k(t), \rho_{k_1}^{r_1,mn}\leq t\leq\rho_{k_1+1}^{r_1,mn}$, $\forall k\in N$, on $\{\omega:\rho^{mn}_k(\omega)\geq\rho_{k_1+1}^{r_1,mn}(\omega)\}$. Let $k\rightarrow +\infty$, it follows that $\theta(\rho^{mn}_k)\rightarrow \theta(\tau^{r_0,mn}_{k_0+1})$, then we have P-a.s.
\begin{equation}\label{eq:3.17}
\sum_{i=1}^2|y_m^i(t)-y_n^i(t)|\leq u^{r_0}_{k_0+1}(t),\quad \rho_{k_1}^{r_1,mn}\leq t\leq\rho_{k_1+1}^{r_1,mn},
\end{equation}
where
\begin{equation*}
u^{r_0}_{k_0+1}(t)=\theta(\tau^{r_0,mn}_{k_0+1})+3\int_t^{\tau^{r_0,mn}_{k_0+1}}(\Phi(u^{r_0}_{k_0+1}(s))+\varepsilon)ds,\\
\end{equation*}
From the arbitrariness of $k_1, r_1$, we can conclude that P-a.s.
\begin{equation}\label{eq:3.18}
\sum_{i=1}^2|y_m^i(t)-y_n^i(t)|\leq u^{r_0}_{k_0+1}(t),\quad \tau_{k_0}^{r_0,mn}\leq t\leq\tau_{k_0+1}^{r_0,mn}.
\end{equation}\\

\textit{Step 5.} By (\ref{eq:3.7}), (\ref{eq:3.18}), similar to the proof of (iv) and (v) in Step 4, we can obtain that for all $r_5>r_0$ and $\forall\tau_{k}^{r_5,mn}\in \{\tau^{r_5,mn}_k\}_{k\in N}$, it still holds that
\begin{equation}\label{eq:3.19}
\sum_{i=1}^2|y_m^i(t,\omega)-y_n^i(t,\omega)|\leq u^{r_5}_k(t,\omega), \quad\tau^{r_0,mn}_{k_0}(\omega)\leq t\leq \tau^{r_0,mn}_{k_0+1}(\omega),
\end{equation}
for P-a.s. $\omega\in\{\omega:\tau_{k_0+1}^{r_0,mn}(\omega)\leq\tau^{r_5,mn}_k(\omega)\}$, where
\begin{equation*}
u^{r_5}_k(t)=\theta(\tau^{r_5,mn}_k)+3\dint_t^{\tau^{r_5,mn}_k}(\Phi(u^{r_5}_k(s))+\varepsilon)ds.
\end{equation*}
Now, we take a subsequence of $\{\tau^{r,mn}_k\}_{r\in Q^+,k\in N}$ denoted by $\{\tau^{mn}_k\}_{k\in N}$ such that $\tau^{mn}_k\rightarrow 1$ as $k\rightarrow +\infty$ and we consider the following equation:
\begin{equation*}
V_k(t)=\theta(\tau^{mn}_k)+3\int_t^{\tau^{mn}_k}(\Phi(V_k(s))+\varepsilon)ds,\\
\end{equation*}
then we have P-a.s. $\sum_{i=1}^2|y_m^i(t)-y_n^i(t)|\leq V_k(t), \tau_{k_0}^{r_0,mn}\leq t\leq\tau_{k_0+1}^{r_0,mn}$, $\forall k\in N$, on $\{\omega:\tau^{mn}_k(\omega)\geq\tau_{k_0+1}^{r_0,mn}(\omega)\}$. Let $k\rightarrow +\infty$, it follows that $\theta(\tau^{mn}_k)\rightarrow \theta(1)=0$ and $P(\{\omega:\tau^{mn}_k(\omega)\geq\tau_{k_0+1}^{r_0,mn}(\omega)\})\rightarrow 1$, then we have P-a.s.
\begin{equation}\label{eq:3.20}
\sum_{i=1}^2|y_m^i(t)-y_n^i(t)|\leq V^{\varepsilon}(t),\quad \tau_{k_0}^{r_0,mn}\leq t\leq\tau_{k_0+1}^{r_0,mn},
\end{equation}
where
\begin{equation*}
V^{\varepsilon}(t)=3\int_t^1(\Phi(V^{\varepsilon}(s))+\varepsilon)ds,\\
\end{equation*}
From the fact that $V^{\varepsilon}(t)\rightarrow 0$ as $\varepsilon\rightarrow 0$ and the arbitrariness of $r_0,k_0$, we can conclude that $(y_n)_{n\geq 0}$ is a Cauchy sequence in $S^{2,d}$ and converges to a process which we denote by $y$. By Step 2, we know that $(z_n)_{n\geq 0}$ is also a Cauchy sequence in $H^{2,d\times m}$ and converges to a process $z$. \\

\textit{Step 6.}
Now, we show the processes $(y,z)$ are the solution to BSDE (\ref{eq:3.3}) on $[0,1]$.  For any $n\geq 1$, we know from (\ref{eq:3.2}) that
\[
y_n(t)=\xi+\int_t^1f_n(s,y_n(s),z_n(s))ds-\int_t^1z_n(s)dB_s, \quad 0\leq t\leq 1.
\]
For a fixed $t$, the sequences $(y_n(t))_{n\geq 0}$ and $(\int_t^1z_n(s)dB_s)_{n\geq 0}$ converge in $L^2(\Omega,{\cal{F}},P)$ towards $y_t$ and $\int_t^1z_sdB_s$, respectively. On the other hand,
\begin{equation*}
\begin{array}{lll}
&&E\left[|\dint_t^1f_n(s,y_n(s),z_n(s))ds-\int_t^1f(s,y(s),z(s))ds|\right]\\
&\leq& E\left[\dint_0^1|f_n(s,y_n(s),z_n(s))-f(s,y_n(s),z_n(s))|ds\right]\\
&&+E\left[\dint_0^1|f(s,y_n(s),z_n(s))-f(s,y(s),z(s))|ds\right].
\end{array}
\end{equation*}
The first term converges to $0$ as $n\rightarrow +\infty$ since $(f_n)_{n\geq 0}$ converges uniformly to $f$. In addition, for any $K\geq 0$, we have
\begin{equation*}
\begin{array}{lll}
&&E\left[\dint_0^1|f(s,y_n(s),z_n(s))-f(s,y(s),z(s))|ds\right]\\
&\leq& E\left[\dint_0^1\Phi(|y_n(s)-y(s)|)1_{\{|y_n(s)-y(s)|\leq K\}}ds\right]\\
&&+E\left[\dint_0^1\Psi(||z_n(s)-z(s)||)1_{\{|z_n(s)-z(s)|\leq K\}}ds\right]\\
&&+E\left[\dint_0^1|f(s,y_n(s),z_n(s))-f(s,y(s),z(s))|1_{\{|y_n(s)-y(s)|\geq K\}}ds\right]\\
&&+E\left[\dint_0^1|f(s,y_n(s),z_n(s))-f(s,y(s),z(s))|1_{\{|z_n(s)-z(s)|\geq K\}}ds\right].
\end{array}
\end{equation*}
Then, after extracting a subsequence, the first term and the second term converge to $0$ as $n\rightarrow +\infty$. Moreover,
\begin{equation*}
\begin{array}{lll}
&&E\left[\dint_0^1|f(s,y_n(s),z_n(s))-f(s,y(s),z(s))|1_{\{|y_n(s)-y(s)|\geq K\}}ds\right]\\
&\leq& \left(E\left[\dint_0^1|f(s,y_n(s),z_n(s))-f(s,y(s),z(s))|^2ds\right]\right)^{\frac{1}{2}}\left(E\left[\dint_0^11_{\{|y_n(s)-y(s)|\geq K\}}ds\right]\right)^{\frac{1}{2}}\\
&\leq&\dfrac{1}{K}\left(E\left[\dint_0^1|f(s,y_n(s),z_n(s))-f(s,y(s),z(s))|^2ds\right]\right)^{\frac{1}{2}}\left(E\left[\dint_0^1|y_n(s)-y(s)|^2ds\right]\right)^{\frac{1}{2}}\\
&\leq& \dfrac{C}{K}\sqrt{E\left[\dint_0^1|y_n(s)-y(s)|^2\right]ds}.
\end{array}
\end{equation*}
The last inequality follows from the convergence of $(y_n)_{n\geq 0}$ $((z_n)_{n\geq 0})$ in $S^{2,d}$ $(H^{2,d\times m})$ and the linear growth of $f$. Therefore, the third term converges to $0$ as $n\rightarrow +\infty$. In the same way, it is easily seen that the forth term also converges to $0$ as $n\rightarrow +\infty$. Consequently, since $y$ is a continuous process, we have
\begin{equation*}
y_t=\xi+\int_t^1f(s,y_s,z_s)ds-\int_t^1z_sdB_s,\quad 0\leq t\leq 1,
\end{equation*}
and we get the existence of the solution of BSDE (\ref{eq:3.3}).\\

\textit{Uniqueness.}
For the uniqueness of the solution to BSDE (\ref{eq:3.3}), the proof is similar to the existence and we omit it.
\hfill $\Box$

\begin{remark}
\label{rmk:3.6}
In Step 3 of the existence, Novikov's condition is satisfied by the fact that $0<\varepsilon_0\leq |z_{mn}(t)|\leq ||z_m(t)-z_n(t)||\leq\frac{1}{\varepsilon_0}$ on  $[{\tau}_{\varepsilon_0,2}^{r_0,mn},{\tau}_{\varepsilon_0,1}^{r_0,mn}]$, hence we can apply Girsanov's theorem on $[0,{\tau}_{\varepsilon_0,1}^{r_0,mn}]$.
\end{remark}

\begin{remark}
\label{rmk:3.7}
For the case of $d>2$, its proof is the same as $d=2$ and the corresponding random differential equation becomes $u(t)=\theta(\tau^{mn})+(d+1)\int_t^{\tau^{mn}}(\Phi(u(s))+\varepsilon)ds$, where $\theta(\tau^{mn})$ is the $(d+1)(\Phi+\varepsilon)$-envelope of the process $\sum_{i=1}^d|y_m^i-y_n^i|$ at stopping time $\tau^{mn}$.
\end{remark}

\begin{remark}
\label{rmk:3.8}
From the proof of Theorem \ref{thm:3.1}, we can see that the condition $\int_{0^+}[\Phi(x)]^{-1}dx= +\infty$ in Assumption 2(i) can be weakened. In fact, if the linear growth function $\Phi$ satisfying $\Phi(0)=0$ and $\Phi(x)>0$ for all $x>0$, moreover, for all $\gamma\geq 0,\varepsilon\geq 0$, the following equation
\[
u^{\gamma,\varepsilon}(t)=\gamma+\int_t^1(\phi(u^{\gamma,\varepsilon}(s))+\varepsilon)ds,\quad 0\leq t\leq 1,
\]
has a unique solution with $\lim\limits_{\varepsilon\rightarrow 0}\lim\limits_{\gamma\rightarrow 0}u^{\gamma,\varepsilon}(t)=0$,
then BSDE (\ref{eq:3.3}) also has a unique solution.
\end{remark}

\section*{Appendix}
We denote
\[
\Phi_k(x):=\sup_{y\in R}\{\Phi(y)-k|x-y|\},\quad k\in N,
\]
from the linear growth of $\Phi$, we know $\Phi_k$ is Lipschitz. Moreover, the sequence $(\Phi_k)_{k\geq 1}$ is non-increasing and converges to $\Phi$. Now , we consider the equation defined recursively as follows:
\begin{equation*}
u^k_{l+1}(t)=E^{mn}\left[\theta({\tau}_{\varepsilon_0,1}^{mn})
+3\dint_t^{{\tau}_{\varepsilon_0,1}^{mn}}(\Phi_k(u^k_l(s))+\varepsilon)ds|{\cal{F}}_t\right],\\
\end{equation*}
with $u^k_1=C, \forall k\geq 1$. From Step 1, we know there exists a nonnegative constant $C$ such that P-a.s. $\sum\limits_{i=1}^2|y_m^i(t)-y_n^i(t)|\leq C, \forall t\in [0,1]$. Then, we can obtain by induction that for all $l,k$,
\begin{equation}\label{eq:3.21}
|y_m^1(t)-y_n^1(t)|+|y_m^2(t)-y_n^2(t)|\leq u^k_l(t), \quad {\tau}_{\varepsilon_0,2}^{mn}\leq t\leq {\tau}_{\varepsilon_0,1}^{mn}.
\end{equation}
Then $\Phi(|y_m(t)-y_n(t)|)\leq \Phi(u^k_l(t))\leq\Phi_k(u^k_l(t))$. In (\ref{eq:3.21}), first, let $l\rightarrow +\infty$ ,then $k\rightarrow +\infty$, we obtain
\[
|y_m^1(t)-y_n^1(t)|+|y_m^2(t)-y_n^2(t)|\leq u^{\varepsilon_0}(t), \quad {\tau}_{\varepsilon_0,2}^{mn}\leq t\leq {\tau}_{\varepsilon_0,1}^{mn}.
\]

\section*{Acknowledgments}
\label{Acknowledgments}
The authors would like to thank Professor Said Hamad\`{e}ne for his helpful comments and suggestions.


\begin{thebibliography}{99}

\bibitem{r1}
\textsc{Dellacherie, C.} and \textsc{Meyer, P.-A.} (1978). \textit{Probabilities and potential}. North-Holland, Amsterdam-New York.
{MR0521810}

\bibitem{r2}
\textsc{El Karoui, N.}, \textsc{Peng, S.} and \textsc{Quenez, M.-C.} (1997). Backward stochastic differential equations in finance. \textit{Math. Finance} \textbf{7} 1-71.
{MR1434407}

\bibitem{r3}
\textsc{Fan, S.}, \textsc{Jiang, L.} and \textsc{Davison, M.} (2010). Uniqueness of solutions for multidimensional BSDEs with uniformly continuous generators. \textit{C.R. Acad. Sci. Paris, Ser, I} \textbf{348} 683-686.
{MR2652498}

\bibitem{r4}
\textsc{Fan, S.}, \textsc{Jiang, L.} and \textsc{Davison, M.} (2013). Existence and uniqueness result for multidimensional BSDEs with generators of Osgood type. \textit{Front. Math. China} \textbf{8} 811-824.
{MR3057221}

\bibitem{r5}
\textsc{Hamad\`{e}ne, S.} (2003). Multidimensional backward stochastic differential equations with uniformly continuous coefficients. \textit{Bernoulli} \textbf{9} 517-534.
{MR1997495}

\bibitem{r6}
\textsc{Kobylanski, M.} (2000). Backward stochastic differential equations and partial differential equations with quadratic growth. \textit{Ann. Probab.} \textbf{28} 259-276.
{MR1782267}

\bibitem{r7}
\textsc{Lepeltier, J.P.} and \textsc{San Martin, J.} (1997). Backward stochastic differential equations with continuous coefficients. \textit{Statist. Probab. Lett.} \textbf{32} 425-430.
{MR1602231}

\bibitem{r8}
\textsc{Lepeltier, J.P.} and \textsc{San Martin, J.} (1998). Existence for BSDE with superlinear quadratic coefficient. \textit{Stoch. Stoch. Rep.} \textbf{63} 227-240.
{MR1658083}

\bibitem{r9}
\textsc{Li, M.} (2015). Backward SDEs with coefficients uniformly continuous in $y$ and linear growth in $z$. Preprint.

\bibitem{r10}
\textsc{Mao, X.} (1995). Adapted solutions of backward stochastic differential equations with non-Lipschitz coefficients. \textit{Stoch. Proc. Appl.} \textbf{58} 281-292.
{MR1348379}

\bibitem{r11}
\textsc{Pardoux, E.} (1999). BSDEs, weak convergence and homogenization of semilinear PDEs. In: Clarke, F., Stern, R. and Sabidussi, G. (Eds), \textit{Nonliner Analysis, Differential Equations and Contral}. NATO Sci. Ser. C Math. Phys. Sci. \textbf{528}, Dordrecht: Kluwer Academic.
{MR1695013}

\bibitem{r12}
\textsc{Pardoux, E.} and \textsc{Peng, S.} (1990). Adapted solution of a backward stochastic differential equation. \textit{Systems Control Lett.} \textbf{14} 55-61.
{MR1037747}

\bibitem{r13}
\textsc{Peng, S.} (1999). Open problems on backward stochastic differential equations. In \textit{Control of distributed parameter and stochastic systems} \textbf{13} 265-273. Kluwer Academic Publ., Boston.
{MR1777419}

\bibitem{r14}
\textsc{Yan, J.} (2004). \textit{Lectures on Measure Theory (Chinese version)}. 2nd ed. Science Press, Bejing.

\end{thebibliography}
\end{document}